\documentclass[10pt]{amsart}
\usepackage[dvips]{graphics}
\usepackage{graphicx}
\usepackage{amsfonts}
\usepackage{euscript}
\usepackage{amssymb}

\usepackage{pstricks}
\usepackage{delarray}
\usepackage{amsmath}
\usepackage{array}
\usepackage{delarray}

\def\R{\mathbb R}
\def\C{\mathbb C}

\def\N{\mathbb N}

\def\({\left(}
\def\){\right)}
\def\[{\left[}
\def\]{\right]}

\def\dem{\noindent  \textbf{Proof. }}

\def\<{\langle}
\def\>{\rangle}

\def\fin{\hfill{$\Box$}}

\def\Ha{\mathcal{H}_{\alpha}}

\def\HA{\operatorname{H_{A(t)}}}
\def\Ha{\operatorname{H_A}}

\def\nabA{\operatorname{\nabla_A}}
\def\tT{\widetilde T}

\def\mg{m\hspace{-0.5mm}g}
\renewcommand\Im{\operatorname{Im}}
\renewcommand\Re{\operatorname{Re}}
\begin{document}

\newtheorem{theorem}{Theorem}
\newtheorem{lemma}{Lemma}[section]
\newtheorem{proposition}[lemma]{Proposition}
\newtheorem{corollary}[lemma]{Corollary}
\newtheorem{definition}[lemma]{Definition}
\newtheorem{remark}[lemma]{Remark}
\newtheorem{example}{ Example}
\newtheorem{assumption}{ Assumption}

\numberwithin{equation}{section}

\title[On non-linear Schr\"odinger equation with magnetic field]
{Remarks on non-linear Schr\"odinger equation with  magnetic fields}
\author[L. Michel]{Laurent Michel}
\address[L. Michel]{Laboratoire J. A. Dieudonn\'e, Universit\'e de
  Nice Sophia-Antipolis, Parc Valrose,
 06108 Nice Cedex 02, France}
\email[L. Michel]{lmichel@math.unice.fr}
\begin{abstract}We study the non-linear Sch\"odinger equation with time depending magnetic field without smallness assumption at infinity. We obtain some results on the Cauchy problem, WKB asymptotics and instability.
\end{abstract}
\maketitle


\section{Introduction}
We consider the non-linear Schr\"odinger equation with magnetic field on $\R^n$
\begin{equation}\label{eq:nlsmagn}
i\partial_tu=\HA u-b^\gamma f(x,u)
\end{equation}
with initial condition
\begin{equation}\label{eq:cond_init_schrod}
u_{|t=t_0}=\varphi.
\end{equation}
Here $$\HA=\sum_{j=1}^n(i\partial_{x_j}-bA_j(t,x))^2,t\in\R,\;\;x\in\R^n$$ is the time-depending Schr\"odinger operator associated to the magnetic potential $A(t,x)=(A_1(t,x),\ldots,A_n(t,x))$, $b\in]0,+\infty[$ is a parameter quantizing the strength of the magnetic field and $\gamma\geq 0$. We sometimes omit the space dependence and write $A(t)$ instead of $A(t,x)$.
The aim of this note is to show that recent improvement in the analysis of non-linear Schr\"odinger equations can be adapted to the case with magnetic field. As an important preliminary, we study the local Cauchy problem for (\ref{eq:nlsmagn}) in energetic space. Let us begin with the general framework of our study.

We suppose that the magnetic potential is a smooth function $A\in C^{\infty}(\R_t\times\R_x^n,\R^n)$ and that it satisfies the following assumption.
\begin{assumption}\label{hyp:dec_magnpot}
\begin{enumerate}
\item $\forall \alpha\in\N^n \displaystyle\sup_{(t,x)\in\R\times\R^n}|\partial_x^\alpha\partial_t A|\leq C_\alpha.$
\item $\forall |\alpha|\geq 1,\displaystyle\sup_{(t,x)\in\R\times\R^n}|\partial_x^\alpha A|\leq C_\alpha.$
\item $\exists\epsilon>0,\forall |\alpha|\geq 1,\displaystyle\sup_{(t,x)\in\R\times\R^n}|\partial_x^\alpha B|\leq C_\alpha\<x\>^{{-1-\epsilon}}$
\end{enumerate}
where $B(t,x)$ is the matrix defined by $B_{jk}=\partial_{x_j}A_k-\partial_{x_k}A_j$.
\end{assumption}
Remark that compactly supported perturbations of linear (with respet to $x$) magnetic potentials satisfy the above hypothesis.

Under Assumption \ref{hyp:dec_magnpot}, the domain $D(\HA)=\{u\in L^2(\R^n_x),\HA u\in L^2(\R^n_x)\}$ does not depend on $t$. Indeed, for $t,t'\in\R$ one has
\begin{equation}\label{eq:depoptps}
\operatorname{H_{A(t')}}=\HA+bW(t,t')(i\nabla_x-bA(t))+b(i\nabla_x-bA(t))W(t,t')+b^2W(t,t')^2
\end{equation}
with $x\mapsto W(t,t',x)=\int_t^{t'}\partial_sA(s,x)ds$ bounded as well as its $x$-derivatives uniformly with respect to $t,t'$ in any compact set.
In fact, the above identity shows that the space
$$H^\beta_{mg}(\R^n)=\{u\in L^2(\R^n),(1+\HA)^{\beta/2} u\in L^2(\R^n)\}$$ does not depend on $t\in\R$. As $D(\HA)=H^2_{mg}(\R^n)$, the above statement is straightforward.
Moreover, the natural norms on this space are equivalent and this equivalence is uniform with respect to the parameter $b$ for close times. More precisely, denoting
$m_A=\sup_{(t,x)\in\R\times\R^n}|\partial_tA(t,x)|$,
we have the following
\begin{proposition}\label{prop:equivnorme}
Suppose that Assumption \ref{hyp:dec_magnpot} is satisfied and let $\beta>0$ and $T>0$. Then,  for all $t,t'\in\R$ such that $|t-t'|\leq b^{-1}T$ and all $u\in H^\beta_{mg}$ we have
$$\|(\operatorname{H_{A(t')}}+1)^\beta u\|_{L^2}\leq (1+2m_A T+m_A^2T^2)^\beta\|(\HA+1)^\beta u\|_{L^2}.$$
\end{proposition}
\dem
It is a straightforward consequence of equation (\ref{eq:depoptps}), Assumption \ref{hyp:dec_magnpot} and the fact that
$(i\nabla_x-bA(t))(\HA+1)^{-1}$ is bounded by $1$ in $L^2$.
\fin

For $\beta\in\N$ we set
\begin{equation}\label{eq:defnormsob}
\|u\|_{H^\beta_{A(t)}}=\|(i\nabla_x-bA(t))^\beta u\|_{L^2}+\|u\|_{L^2}.
\end{equation}
This norm is clearly equivalent (uniformly with respect to $b$) to
$\|(1+\operatorname{H_{A(t)}})^{\beta/2}u\|_{L^2}$.
In regard of Proposition \ref{prop:equivnorme} we define the magnetic Sobolev norm by
$$\|u\|_{H^\beta_{mg}}=\|u\|_{H^\beta_{A(t_0)}}.$$

Under Assumption \ref{hyp:dec_magnpot} it is well-known (see \cite{Pa83}, Th 4.6, p143 or \cite{Ya91}) that
for $\varphi\in H^1_{mg}$, the linear Schr\"odinger equation
\begin{equation}\label{eq:lsmagn}
i\partial_tu=\HA u,\;\;u_{|t=s}=\varphi
\end{equation}
has a solution $U_0(t,s)\varphi$. The operator $U_0(t,s)$ maps $H^1_{mg}$ into itself, is continuous from $L^2$ into $L^2$ and from $H^1_{mg}$ into $H^1_{mg}$. Moreover, $U_0(t,s)\varphi$ is the unique $H^1_{mg}$ valued solution of (\ref{eq:lsmagn}) and $U_0(t,s)$ is unitary.

The first aim of this paper is to solve the Cauchy problem for the non-linear equation in the most appropriate space.
We state the
assumptions on the non-linearity $f$. We suppose that
$f:\R^n\times\C\rightarrow \C$ is a measurable function such that
\begin{assumption}\label{hyp:nonlinearite}
\begin{enumerate}
\item $f(x,0)=0$ almost every where.\\
\item $\exists M\geq 0,\alpha\in[0,\frac 4 {n-2}[$ ($\alpha\in[0,\infty[$ if $n=1,2$) such that
$$|f(x,z_1)-f(x,z_2)[\leq M(1+|z_1|^\alpha+|z_2|^\alpha)|z_1-z_2|$$
for almost all $x\in\R^n$ and all $z_1,z_2\in\C$.\\
\item $\forall z\in\C, f(x,z)=(z/|z|)f(x,|z|)$
\end{enumerate}
\end{assumption}
Remark that these assumptions are often used in the case $A=0$. More precisely, in the case $A=0$, the second property of the above assumption corresponds to a subcritical non-linearity with respect to $H^1$.

Let us introduce some energy functional associated to these non-linerarities. We define 
\begin{equation*}
\begin{split}
 &F(x,z)=\int_0^{|z|}f(x,s)ds,\;G(u)=\int_{\R^n}F(x,u(x))dx
\end{split}
\end{equation*}
and for $t\in\R$ and $u\in H^1_{mg}$ we define the energy
\begin{equation*}
 \begin{split}
  &E(b,t,u)=\int_{\R^n}\frac 1 2 |(i\nabla_x-bA(t,x))u(x)|^2dx-b^{\gamma}G(u).
 \end{split}
\end{equation*}
Formally, it is not hard to see that any sufficiently regular solution of (\ref{eq:nlsmagn}), (\ref{eq:cond_init_schrod}), enjoys the following energy evolution law:
$$E(b,t,u)=E(b,0,\varphi)-\Re\int_0^t\<\partial_sA(s)u(x),(i\nabla-A(s))u(s)\>_{L^2}ds.$$
Therefore, the natural space to solve (\ref{eq:nlsmagn}), (\ref{eq:cond_init_schrod}) seems to be $H^1_{mg}$.

Now we are in position to state our first result.
\begin{theorem}\label{th:pbcauchy_nlsmagn}
Suppose that Assumptions \ref{hyp:dec_magnpot} and
\ref{hyp:nonlinearite} are satisfied and let $\varphi\in
H^1_{mg}$. Then, there exists $T_b,T^b>0$ and a unique $u\in
C(]-T_b,T^b[,H^1_{mg})\cap C^1(]-T_b,T^b[,H^{-1}_{mg})$ solution of
(\ref{eq:nlsmagn}). Moreover, either $T_b=\infty$ (resp. $T^b=\infty$), or
$\lim_{t\rightarrow -T_b}\|u(t)\|_{H^1_{\mg}}=\infty$ (resp. $\lim_{t\rightarrow T^b}\|u(t)\|_{H^1_{\mg}}=\infty$) and
\begin{equation}\label{eq:consmasse}
\|u(t)\|_{L^2}=\|\varphi\|_{L^2},
\end{equation}
\begin{equation}\label{eq:loienerg}
E(b,t,u)=E(b,0,\varphi)-\Re\int_0^t\<\partial_sA(s)u(x),(i\nabla-A(s))u(s)\>_{L^2}ds,
\end{equation}
for all $t\in ]-T_b,T^b[$.
Additionally,
there exists $\epsilon>0$ such that, for all $b>0$ and $\varphi\in
H^1_{mg}$ such that $\|\varphi\|_{H^1_{mg}}\leq C b$, we have
$T_b,T^b\geq \epsilon b^{-\delta}$ with
$\delta=\max(1,2\gamma,\frac{2\gamma}{\alpha})$.
\end{theorem}

Let us make a few remarks on this result. The Cauchy problem for non-linear Schr\"odinger equation has a long story. In absence of magnetic field there are numerous results; see for instance \cite{GiVe79,GiVe85}, \cite{CaWe88}. 

In presence of magnetic field, the behavior of $A$ when $|x|$ becomes large plays an important role. In the case where the magnetic potential $A$ is bounded, the spaces $H^1_{mg}$ and $H^1$ coincide and the Cauchy problem can be solved in $H^1$ using usual techniques. If the magnetic field is unbounded , it is not possible to solve the Cauchy problem in $H^1$ as multiplication by $A$ is not bounded on $L^2$.

To avoid this difficulty some authors work in the weighted Sobolev space $\Sigma=\{u\in H^1(\R^n),\;(1+|x|)u\in L^2,\}$ (see for instance \cite{Bo91}, \cite{NaSh05}). In particular, they require some decay of the initial data at infinity.

 In the case of \cite{Bo91}, this decay is required because the author use dispersive properties for the Laplacian instead of $\HA$. In \cite{NaSh05} the author use magnetic Strichartz estimates but their method based on fixed-point theorem is not adapted to the magnetic context and requires decay of the solution at infinity.

On the other hand, there exists also of a result of Cazenave and Esteban \cite{CaEs88} dealing with the special case where the magnetic field $B$ is constant (and hence, $A$ is linear with respect to $x$).
In a way, this paper is more satisfactory as they need only $u_0$ to belong to the energy space. Nevertheless, their result applies only to constant magnetic field.

Our theorem is, then a generalization of the above results. Before going further, let us remark that for unbounded $A$, the spaces $H^1$, $H^1_{mg}$ and $\Sigma$ are different. First, it is evident that $\Sigma$ is contained in $H^1\cap H^1_{mg}$. Let us give an example where $\Sigma$ is strictly contained in $H^1_{mg}$. For this purpose, we restrict ourseleves to the case where the dimension $n=2$ and consider the magnetic potential $A(x,y)=(y,x)$. Let $g\in H^1(\R^2)$ be such that $|x|g\notin L^2$, then a simple calculus shows that $f(x,y)=g(x,y)e^{-ixy}$ belongs to $H^1_{mg}\setminus \Sigma$.

In the case of defocusing non-linearities the energy law implies the following result.

\begin{corollary}
Suppose that $F(x,z)\leq 0$ for all $x,z$, then $T_b,T^b=+\infty$.
\end{corollary}
\dem
For $F\leq0$, we deduce from (\ref{eq:loienerg}) and Cauchy-Schwarz inequality, that 
$$\|(i\nabla-A(t))u(t)\|_{L^2}\leq C_1+C_2\int_0^t\|(i\nabla-A(s))u(s)\|_{L^2}ds$$
for some fixed constant $C_1,C_2>0$. Hence, Gronwall Lemma shows that $\|(i\nabla-A(t))u(t)\|_{L^2}$ remains bounded on any bounded time-interval. Using (\ref{eq:consmasse}) and the characterization of $T_b$, we obtain the result.
\fin

The next section contains the proof of Theorem \ref{th:pbcauchy_nlsmagn}. In section 3 we give some qualitative results on the solution of (\ref{eq:nlsmagn}) in the limit $b\rightarrow \infty$. More precisely, we can construct WKB solutions and prove instability results with respect to initial data and parameter $b$.

\section{Cauchy problem in the energy space}

The proof of theorem \ref{th:pbcauchy_nlsmagn} relies on the Strichartz estimates proved in \cite{Ya91} for the problem
\begin{equation}\label{eq:lsmagn_inhom}
 i\partial_tu=\HA u+g(t),\;\; u_{|t=s}=\varphi
\end{equation}

\begin{theorem}\label{th:strichmagn}{\bf(Yajima)} Let $I$ be a finite real interval, $(q,r)$ and $(\gamma_j,\rho_j), j=1,2$ be such that
$r,\rho_j\in[2,\frac {2n}{n-2}[$, $\frac 2 q= n(\frac 1 2 -\frac 1 r)$ and $\frac 2 {\gamma_j}= n(\frac 1 2 -\frac 1 {\rho_j})$. Let
$g_j\in L^{\gamma'_j}(I,L^{\rho'_j}(\R^n_x)),j=1,2$, where $\gamma'_j,\rho'_j$ are the conjugate exposant of $\gamma_j, \rho_j$. Then the solution $u$ to
(\ref{eq:lsmagn_inhom}) with $g=g_1+g_2$ satisfies
\begin{equation}\label{eq:strichmagn}
 \|u\|_{L^q(I,L^r(\R^n_x))}\leq C (\|g_1\|_{L^{\gamma'_1}(I,L^{\rho'_1}(\R^n_x))}+\|g_2\|_{L^{\gamma'_2}(I,L^{\rho'_2}(\R^n_x))}+\|\varphi\|_{L^2(\R^n)})
\end{equation}
where the constant $C$ depends only on the length of $I$ and the
constant $C_\alpha$ of Assumption \ref{hyp:dec_magnpot}.
\end{theorem}
\dem In the case $g=0$ it is exactly Theorem 1 of \cite{Ya91}. In
the general case it suffices to work as in the proof of Proposition
2.15 of \cite{BuGeTz04} using a celebrated result of Christ and
Kiselev \cite{ChKi01}. The fact that the constant $C$ depends only
on the $C_\alpha$ is a direct consequence of the construction of
Yajima \cite{Ya91}. \fin

\begin{remark}
In the case where the magnetic potential is not regular, there are some recent results of A. Stefanov \cite{St05} and Georgiev-Tarulli \cite{GeTa06} which provide Strichartz estimates under smallness assumption on the magnetic fields. This should lead to the corresponding existence and uniqueness result for NLS in the case of small magnetic field. This could also have consequences on the well-posedness of the Schr\"odinger-Maxwell system (see \cite{GiVe06}, \cite{NaWa05}, \cite{Ts93} for results on this topics).
\end{remark}

It is important to notice that Theorem \ref{th:pbcauchy_nlsmagn} is
not a straightforward consequence of the above Strichartz estimate.
Indeed, if we try to apply a fixed point method to equation
(\ref{eq:nlsmagn}), a problem occurs when we try to control the norm
of the non-linearity in the $H^1_{mg}$ norm. Consider for instance
the case $f(u)=|u|^2u$, then
$$(i\nabla_x-bA(t))(|u|^2u)=|u|^2(i\nabla_x-bA(t))(u)+ui\nabla_x(|u^2|).$$
The first therm of the right hand side of this equality will be controlled by $\|u\|_{H^1_{mg}}$, whereas in the second term, as $A(t,x)$ is not bounded with respect to $x$, there is no chance to control $i\nabla_x(|u^2|)$ by $(i\nabla_x-bA(t))(|u^2|)$. 
For the same reason it does not seem easy to solve the Cauchy problem in magnetic Sobolev spaces of high degree.

To overcome this difficulty, we work as in \cite{CaWe88}, \cite{CaEs88} and approximate the solution of (\ref{eq:nlsmagn}) by solution of a non-linear Schr\"odinger equation with non-linearity linearized  at infinity. In the work of Cazenave and Weissler, the main tool to justify the approximation is an energy conservation. In our case, the Hamiltonian depends on time, so that the energy is not conserved. Nevertheless, the error term is controled by the $H^1_{mg}$-norm so that it is possible to implement the same strategy. 
Another difference involved by the dependance with respect to time of the Hamiltonian is that usual techniques to solve the Cauchy problem with regular initial data and nice non-linearities can not apply in our context. Therefore, additionnaly to the approximation of the non-linearity, we have to introduce an approximation of the magnetic field itself and justify the convergence to our initial problem.
\\

Let us introduce the approximated nonlinearities used in the sequel.
Following \cite{CaWe88}, we decompose $f=\tilde f_1+\tilde f_2$
 with
\begin{equation}\label{eq:def_tf1}
\tilde f_1(x,z)=1_{\{|z|\leq1\}}f(x,z)+1_{\{|z|\geq1\}}f(x,1)z
\end{equation}
and
\begin{equation}\label{eq:def_tf2}
\tilde f_2(x,z)=1_{\{|z|\geq1\}}(f(x,z)-f(x,1)z).
\end{equation}
Next we define $f_m=\tilde f_1+\tilde f_{2,m}$ where
\begin{equation}\label{eq:def_tf2m}
\tilde f_{2,m}(x,z)=1_{\{|z|\leq m\}}\tilde f_2(x,z)+1_{\{|z|\geq
m\}}\tilde f_2(x,m)\frac z m
\end{equation}
 Remark that these
functions satisfy Assumption \ref{hyp:nonlinearite}.
We consider also the energy functional associated to these approximated non-linearities.
We define
\begin{equation}
\begin{split}
F_m(x,z)=\int_0^{|z|}f_m(x,s)ds,\;G_m(u)\int_{\R^n}F_m(x,u(x))dx
\end{split}
\end{equation}
and for $t\in\R$ and $u\in H^1_{mg}$ we set
\begin{equation}
 \begin{split}
E_m(b,t,u)=\int_{\R^n}\frac 1 2 |(i\nabla_x-bA(t,x))u(x)|^2dx-G_m(u).
 \end{split}
\end{equation}

Finally, we remark that replacing
the magnetic potential $A(t,x)$ by $A(t+t_0,x)$ it suffices to prove
Theorem \ref{th:pbcauchy_nlsmagn} for $t_0=0$.

On the other hand, to enlight the notations we prove the theorem in the particular case $b=1$. To get the general case it suffices to keep track of $b$ along the proof. We will also restrict our study to $t\geq 0$, the other case being treated by reversing time in the equation.

\subsection{Preliminary results}
In the sequel, we will need Sobolev embeddings in the magnetic context. In this subsection, $A$ is a magnetic potential satisfying Assumption \ref{hyp:dec_magnpot}.
\begin{lemma}\label{lem:injecsobmagn}
Let $0<s<\frac n 2$ and $p_s=\frac{2n}{n-2s}$, then $H^s_{A}$ is
continuously embedded in $L^{p}(\R^n)$ for all $p\in [2,p_s]$ and
there exists $C>0$ independent of $A$ such that
\begin{equation}\label{eq:inegsob}
\|u\|_{L^p}\leq C\|u\|_{H^s_{A}}
\end{equation}

\end{lemma}
\dem From the diamagnetic inequality (see \cite{AvHeSi78}), we know that
almost everywhere we have
$$|u|=|(\Ha+1)^{-\frac s 2}(\Ha+1)^{\frac s 2}u|\leq
(-\Delta+1)^{-\frac s 2}|(\Ha+1)^{\frac s 2}u|.$$ Taking the $L^p$
norm, the result follows from standard Sobolev inequalities.
 \fin

Next we prove a technical result on the non-linearity. 

\begin{proposition}\label{prop:estim_nonlin}
 Let $M>0$, $r_1=\rho_1=2$ and $r_2=\rho_2=\alpha+2$  then
\begin{enumerate}

\item the sequence $(\tilde f_{2,m}(.,u))_{m\in\N^*}$ converges to $\tilde f_2(.,u)$ in $L^{\rho_2'}(\R^n)$ uniformly with respect to $u\in H^1_{A}$ such that
$\|u\|_{H^1_{A}}\leq M$.
\item there exists $C(M)>0$ independent of $A$ such that for all $m\in\N^*$ and for all $u,v\in H^1_{A}$ with
$\max(\|u\|_{H^1_{A}},\|v\|_{H^1_{A}})\leq M$ we have
$$\|\tilde f_1(.,u)-\tilde f_1(.,v)
\|_{L^{\rho'_1}(\R^n)}\leq C(M) \|u-v\|_{L^{r_1}}$$
\begin{equation*}
\begin{split}
 \|\tilde f_{2,m}(.,u)-\tilde f_{2,m}(.,v)\|_{L^{\rho'_2}(\R^n)}+\|\tilde f_2(.,u)-\tilde f_2(.,v)
\|_{L^{\rho'_2}(\R^n)}\leq C(M)\|u-v\|_{L^{r_2}(\R^n)}
\end{split}
\end{equation*}

\end{enumerate}
\end{proposition}

\dem We follow the method of Example 3 in \cite{CaWe88}. Taking
$\chi$ the characteristic function of the set
$\{x\in\R^n\,|\;|u(x)|>m\}$ and using Assumption
\ref{hyp:nonlinearite}, we have
\begin{equation}\label{eq:appnl1}
\|\tilde f_2(u)-\tilde f_{2,m}(u)\|_{L^{\rho'_2}(\R^n)}\leq 2 \|\chi
|u|^{\alpha+1}\|_{L^{\rho_2'}}=2\|\chi u\|_{L^{\alpha+2}}
^{\alpha+1}.
\end{equation}
On the other hand, using H\"older inequality and Lemma
\ref{lem:injecsobmagn} we get for $p=\frac{2n}{n-2}$,
\begin{equation}\label{eq:appnl2}
\|u\|_{H^1_{A}}\geq C \|\chi u\|_{L^p}\geq Cm^{1-\frac{\alpha
+2}p}\|\chi u\|_{L^{\alpha +2}}^{\frac {\alpha +2}p}.
\end{equation}
As $\alpha<\frac 4 {n-2}$ then $1-\frac{\alpha}{p+2}>0$.
Combining equations (\ref{eq:appnl1}) and (\ref{eq:appnl1}), we
obtain the first point of the proposition.

The second assertion follows, as in example 3 in \cite{CaWe88}, from H\"older inequality,
Assumption \ref{hyp:nonlinearite} and Lemma \ref{lem:injecsobmagn}. The fact that the constant $C(M)$ is independent of the magnetic fields follows from the uniformity of the constant in Lemma \ref{lem:injecsobmagn}.
\fin

\begin{lemma}\label{lem:estim_energnl} Let $T>0$ and $\gamma_k,\;k=1,2$ be
defined by $\frac 2 {\gamma_k}=n(\frac 1 2-\frac 1{\rho_k})$. For
$M>0$ there exists a constant $C(M)$ independent of $A$, such that
for all $u,v\in H^1_{A}$ with $\|u\|_{H^1_{A}}\leq M$ and
$\|v\|_{H^1_{A}}\leq M$ we have
$$|G(u)-G(v)|+|G_m(u)-G_m(v)|\leq C(M)(\|v-u\|_{L^2}+\|v-u\|_{L^2}^{\nu}),$$
with $\frac 2 \nu= \frac n 2 -\frac n {\alpha+2}$ and for all
$u,v\in L^\infty([0,T]H^1_{A})$,
$$\|\tilde f_1(.,u)-\tilde f_1(.,v)
\|_{L^{\gamma'_1}([0,T],L^{\rho'_1}(\R^n))}\leq C(M)T
\|u-v\|_{L^{\gamma_1}([0,T],L^{r_1}(\R^n))}.$$
\begin{equation*}
\begin{split}
 \|\tilde f_{2,m}(.,u)-\tilde f_{2,m}(.,v)\|_{L^{\gamma'_2}([0,T],L^{\rho'_2}(\R^n))}&+\|\tilde f_2(.,u)-\tilde f_2(.,v)
\|_{L^{\gamma'_2}([0,T],L^{\rho'_2}(\R^n))}\\
&\leq
C(M)T^{\frac{\gamma_2-1}{\gamma_2}}\|u-v\|_{L^{\gamma_2}([0,T],L^{r_2}(\R^n))}
\end{split}
\end{equation*}

 Moreover, $G_m\rightarrow G$ as
$m\rightarrow\infty$ uniformly on bounded sets of $H^1_{A}$.
\end{lemma}
\dem
Remark that $G(u)=\int_0^1\<f(x,su),u\>_{L^2}ds$ and $G_m(u)=\int_0^1\<f_m(x,su),u\>_{L^2}ds$ and copy the proof of Lemma 3.3 in \cite{CaWe88}, replacing classical Sobolev inequalities by Lemma  \ref{lem:injecsobmagn} and using Proposition \ref{prop:estim_nonlin}.
\fin
\vspace{0.5cm}

We are now in position to prove the uniqueness part of Theorem \ref{th:pbcauchy_nlsmagn}.
\begin{proposition}\label{prop:unicite}
Let $T>0$ and $u,v\in C([0,T[,H^1_{mg})\cap C^1([0,T[,H^{-1}_{mg})$ be solution of (\ref{eq:nlsmagn}). Then $u=v$.
\end{proposition}

\dem Let $u,v\in C([0,T[,H^1_{mg})\cap C^1([0,T[,H^{-1}_{mg})$ be solution of (\ref{eq:nlsmagn}), and set $w=v-u$. Then $w(0)=0$ and
$$i\partial_t w-\HA w=\widetilde f_1(u)-\widetilde f_1(v)+\widetilde f_2(u)-\widetilde f_2(v).$$
Let $r\in[2,\frac 2 {n-2}]$ and $q>2$ such that $\frac 2 q=n(\frac 1 2-\frac 1 r)$. Apply Theorem \ref{th:strichmagn} together with Lemma \ref{lem:estim_energnl}, we get
$$\|w\|_{L^q([0,T[,L^r)}\leq C(T+T^{\gamma_2})(\|w\|_{L^\infty([0,T[,L^2)}+\|w\|_{L^{\gamma_2}([0,T[,L^{\rho_2})})$$
where $\gamma_2=\frac{\alpha+1}{\alpha+2}$ and $\frac 2 {\gamma_2}=n(\frac 1 2-\frac 1 {\rho_2})$. As we can alternatively take $(q,r)$ to be equal to $(2,\infty)$ and $(\gamma_2,\rho_2)$, we get the announced result by summing the obtained inequalities and making $T>0$ small enough.
\fin
\subsection{Autonomous case}
In this section we explain briefly how to solve the Cauchy problem in $H^1_{mg}$ when the magnetic field $A(t,x)=A(x)$ is time independent.
In this context, the functional $E$ does not depend  on time and formally we have the following conservation of energy. Suppose that $u$ is solution of (\ref{eq:nlsmagn}) then
\begin{equation*}
E(b,u(t))=E(b,\varphi),\;\forall t.
\end{equation*}
More precisely, we prove the following
\begin{proposition}\label{prop:pbcauchy_nlsmagn_auton}
Let $M>0$ and $C_\alpha,\;\alpha\in\N^n$ a family of finite positive
numbers. There exists $T>0$ depending only on $M$ and the $C_\alpha$
such that for all $A$ satisfying $\partial_tA=0$ and Assumption \ref{hyp:dec_magnpot}
with $C_\alpha$ and for all $\varphi\in H^1_{A}$ such that
$\|\varphi\|_{H^1_A}\leq M$, there exists a unique $u\in
C^0([0,T[,H^1_{A})\cap C^1(]0,T[,H^{-1}_{A})$ maximal solution of
\begin{equation*}
i\partial_t u=\Ha u+f(x,u)
\end{equation*}
with initial condition $u_{|t=0}=\varphi$.
Moreover, for all $t\in[0,T[$ we have
\begin{equation*}
E(b,u(t))=E(b,\varphi),
\end{equation*}
and if $T<\infty$ then $\lim_{t\rightarrow T}\|u\|_{H^1_A}=\infty$.
\end{proposition}

The proof is slight adaption of \cite{CaWe88}, \cite{CaEs88} to our
context. We need also to investigate the dependence of the existence
time with respect to the magnetic field. However, the scheme of
proof is the same and consists to consider an approximate problem
and justify convergence on fixed time intervals. Let us give the
main steps of the proof.

{\bf Step 1.} Let $f_m$ be defined by (\ref{eq:def_tf1}), (\ref{eq:def_tf2}), (\ref{eq:def_tf2m}) and let $A$ be a magnetic field satisfying the above hypotheses. Consider the problem
\begin{equation}\label{eq:nlsmagn_app_aut}
i\partial_t u=\Ha u+f_m(x,u),\;\;u_{t=0}=\varphi
\end{equation}
with $\varphi\in H^1_A$. We have the following
\begin{lemma}\label{lem:exist_app_aut}
 Let $\varphi\in H^1_{A}$, then there exists $\tau_{m,A}>0$ such that  there exists $u_m\in C([0 ,\tau_{m,A}[,H^1_{A})\cap C^1([0,\tau_{m,A}[,H^{-1}_{A})$
solution of (\ref{eq:nlsmagn_app_aut}). Moreover we have for all $t\in[0,\tau_{m,A}[$,
\begin{equation}
E_m(u_{m})=E_m(\varphi)
\end{equation}
and
\begin{equation}
 \|u_m(t)\|_{L^2}=\|\varphi\|_{L^2}.
\end{equation}
\end{lemma}
\dem
The proof is the same as that of Lemma 3.5 of \cite{CaWe88}, replacing usual derivatives by magnetic derivatives.
\fin

{\bf Step 2.} We show that the existence time $\tau_{m,A}$ can be bounded from below uniformly with respect to $m\in\N$ and $A$ satisfying Assumptions of the above proposition.
\begin{lemma}\label{lem:minortps_aut}
 Let $M>0$. There exists $T_1=T_1(M)>0$ such that for all $m\in\N$, all $A$ satisfying Assumption \ref{hyp:dec_magnpot}
 and all $\varphi\in H^1_A$ with $\|\varphi\|_{H^1_A}\leq M$ we have
$$\|u_m\|_{L^\infty([0,T_1],H^1_A)}\leq 2\|\varphi\|_{H^1_A}.$$
\end{lemma}
\dem The proof is exactly the same as in Lemma 3.6 of \cite{CaWe88},
making use of  Lemma
\ref{lem:exist_app_aut} (in particular, we use strongly the conservation of energy) and Proposition \ref{prop:estim_nonlin} to get uniformity with respect to $A$. \fin

{\bf Step 3.} The final step is to prove convergence of the $u_m$ to
solution of the initial problem. First we prove convergence in $L^2$.
\begin{lemma}\label{lem:cauchyL2aut}
Let $M>0$ and $C_\alpha,\;\alpha\in\N^n$ a family of finite positive
numbers. There exists $T_2>0$ depending only on $M$ and the
$C_\alpha$ such that for all $A$ satisfying Assumption
\ref{hyp:dec_magnpot} with $C_\alpha$ and for all $\varphi\in
H^1_{A}$ such that $\|\varphi\|_{H^1_A}\leq M$, such that
$(u_m)_{m\in\N}$ is a Cauchy sequence in $C([0,T_2],L^2)$.
\end{lemma}
\dem The proof is the same as in \cite{CaWe88}, making use of
Theorem \ref{th:strichmagn}, Lemma \ref{lem:estim_energnl},
Proposition \ref{prop:estim_nonlin} and Lemma \ref{lem:exist_app_aut}.
\fin

Now, we can complete the proof of Theorem \ref{th:pbcauchy_nlsmagn}.
 We denote $u$ the limit of $u_m$ in $C([0,T_2],L^2)$.
From Lemma \ref{lem:minortps_aut}, it follows that $u\in
L^{\infty}([0,T_2],H^1_{A})$ and by Lemma \ref{lem:injecsobmagn},
$u_m$ converges to $u$ in $C([0,T_2],L^r)$ for all $r\geq 2n/(n-2)$.
Hence, it follows from Proposition \ref{prop:estim_nonlin} that
$f_m(u_m)$ converges to $f(u)$ in $C([0,T_2],H^{-1}_{A})$ and $u$
solves \ref{eq:nlsmagn} in $L^{\infty}([0,T_2],H^{-1}_{A})$.
Moreover, combining Lemma \ref{lem:estim_energnl} and
\ref{lem:exist_app_aut} we prove that
$$E(b,t,u)=E(b,0,\varphi).$$
This shows that $u\in C([0,T_2],H^1_{mg})$ and hence $u\in
C^1([0,T_2],H^{-1}_{A})$.

\subsection{Cauchy problem in the time-depending case}
We suppose now that $A(t,x)$ satisfies Assumption
\ref{hyp:dec_magnpot}. The strategy of proof is the same as in autonomous case and we first consider the problem
\begin{equation}\label{eq:nlsmagn_app}
i\partial_t u=\HA u+f_m(x,u),\;\;u_{t=0}=\varphi
\end{equation}
At least formally, we can see that the energy of the solution of this equation satisfies the following rule
\begin{equation}\label{eq:energiecons}
E(t,u)=E(0,\varphi)-\Re\int_0^t\<\partial_sA(s)u(s),(i\nabla_x-A(s))u(s)\>ds.
\end{equation}
This will replace the energy conservation in our approach.
On the other hand another problem occurs if we try to apply the proof of \cite{CaWe88}. Indeed, the first step should be to obtain a generalization of Lemma \ref{lem:exist_app_aut} in the time depending framework. Following the proof of Lemma 3.5 in \cite{CaWe88}, we should regularize the initial data and solve the Cauchy problem in $H^2_{mg}$. The issue is that contrary to the autonomous case, the existence of smooth solution is not easy to prove. Indeed, the key point in the approach of \cite{CaWe88} is that for any $g\in C([0,T],H^1)$ Lipschitz continuous with respect to time, the function $v(t)=\int_0^tU_0(t,s)g(s)ds$ is also Lipschitz continuous with respect to time. Such a result is easely proved in the autonomous case as the identity $U_0(t+h,s)=U_0(t,s-h)$ permits to use the assumption on $g$. This fails to be true in the time-depending case. For this reason, we prove the existence in $H^1_{mg}$ in a direct way.

\subsubsection{Existence of solution for approximated problem}
In the case where the magnetic potential depends on time, we can not use the method of \cite{CaWe88} to prove existence of solution on (\ref{eq:nlsmagn_app}) in $H^1_{mg}$.  However we can prove the following.

\begin{proposition}\label{prop:exist_app}
 Let $\varphi\in H^1_{mg}$, then there exists $\tT>0$ such that  there exists
 $u_m\in C([0,\tT[,H^1_{mg})\cap C^1([0,\tT[,H^{-1}_{mg})$
solution of (\ref{eq:nlsmagn_app}). Moreover we have for all
$t\in[0,\tT]$,
\begin{equation}\label{eq:estimenergie_app}
E_m(t,u_{m})=E_m(t,\varphi)-\Re\int_0^t\<\partial_sA(s)u(s),(i\nabla_x-A(s))u(s)\>ds.
\end{equation}
and
\begin{equation}\label{eq:conservmasse}
 \|u_m(t)\|_{L^2}=\|\varphi\|_{L^2}.
\end{equation}
\end{proposition}
\dem
The method consists in approximating the magnetic potential $A(t,x)$ by potentials which are
piecewise constant with respect to time.
More precisely, remark that thanks to Assumption \ref{hyp:dec_magnpot} and
Proposition \ref{prop:pbcauchy_nlsmagn_auton} there exists $T_2=T_2(M)>0$ such that for all $t_0\in[0,T_2]$ the Cauchy problem
$$i\partial_tu=\operatorname{H_{A(t_0)}}u(t)+f_m(u(t)),\;u_{|t=t_0}=\varphi$$ can be solved in
$C([t_0,t_0+T_2],H^1_{A(t_0)})$ for all initial data such that $\|\varphi\|_{H^1_{A(t_0)}}\leq M$.

Let $T\in]0,T_2[$ and for $n\in\N^*,\;k\in\{0,\ldots,n-1\}$ define $t_n^k=\frac {kT}n$. We set
$A_n(t,x)=A(t_n^k,x),\;\forall t\in[t_n^k,t_n^{k+1}[$ and $A_n(T,x)=A(T,x)$. Next, we define the Hamiltonian
$H_n=(i\nabla_x-A_n)^2$ and we look for solutions $u_{n,m}$ of
\begin{equation}\label{eq:nlsmagn_loct}
i\partial_tu=H_nu+f_m(u),\;u_{|t=0}=\varphi.
\end{equation}
From uniqueness in the autonomous case, such a function is given by
\begin{equation}\label{eq:def_unm}
u_{n,m}(t,x)=\sum_{k=0}^{n-1}1_{[t_n^k,t_n^{k+1}[}(t)v_{k,n,m}(t,x)
\end{equation}
where $v_{k,n,m}(t,x)$ is defined as follows.
We choose $v_{0,n,m}$ to be solution of
\begin{equation}\label{eq:v0nm}
\left\{\begin{array}{cc}
i\partial_tv_{0,n,m}=(i\nabla_x-A(t_n^0,x))^2v_{0,n,m}+f_m(v_{0,n,m})&\\
v_{0,n,m}(t_n^0,x)=\varphi(x)\phantom{*****************}
\end{array}\right.
\end{equation}
and for $k\geq 1$, $v_{k,n,m}(t,x)$ is the solution of
\begin{equation}\label{eq:vknm}
\left\{\begin{array}{cc}
i\partial_tv_{k,n,m}=(i\nabla_x-A(t_n^k,x))^2v_{k,n,m}+f_m(v_{k,n,m})&\\
v_{k,n,m}(t_n^k,x)=v_{k-1,n,m}(t_n^k,x).\phantom{*************}
\end{array}\right.
\end{equation}
Thanks to Proposition \ref{prop:pbcauchy_nlsmagn_auton}, the
function $v_{k,n,m}$ are well defined and belong to
$C^0([t_n^k,t_n^k+T_2],H^1_{mg})$ and satisfy the following
conservation equations
\begin{equation*}
E_{n,m}(t,v_{k,n,m}(t))=E_{n,m}(t_n^k,v_{k,n,m}(t_n^k))
\end{equation*}
for all $k=0,\ldots,n-1$, $t\in[t_n^k,t_n^{k+1}[$ and where for all $w\in H^1_{mg}(\R^n)$,
$$E_{n,m}(t,w)=\frac 1 2\int_{\R^n}|(i\nabla_x-A_n(t,x))w(x)|^2dx-G_m(w).$$
Let us write $A(t_n^k,x)=A(t_n^{k-1},x)+W_{n,k}(x)$ with
$W_{n,k}(x)=\frac 1 2\int_{t_n^{k-1}}^{t_n^k}\partial_tA(t,x)dt$ and
use $v_{k,n,m}(t_n^k,x)=v_{k-1,n,m}(t_n^k,x)$, then
\begin{equation}\label{eq:cons_en_disc}
\begin{split}
E_{n,m}(t_n^k,u_{n,m}(t_n^k))&=E_{n,m}(t_n^{k-1},u_{n,m}(t_n^{k-1}))\\
&-\int_{t_n^{k-1}}^{t_n^k}\Re\<(i\nabla_x-A(t_n^{k-1}))u_{n,m}(t_n^{k-1}),\partial_tA(t,x)u_{n,m}(t_n^{k-1})\>dt\\
&+\|W_{n,k}u_{n,m}(t_n^{k-1})\|^2_{L^2}.
\end{split}
\end{equation}
Thanks to Assumption \ref{hyp:dec_magnpot} and conservation of mass,
we have
$\|W_{n,k}u_{n,m}(t_n^{k})\|^2_{L^2}=O(\frac{\|\varphi\|^2_{L^2}}
{n^2})$ uniformly with respect to $k,n,m$.

Hence, taking the sum of equations (\ref{eq:cons_en_disc}) for
$k=1,\ldots,k_0$ with $k_0=[\frac {nt}T]$, and using the fact that
the energy is constant on $[t_n^{k_0},t_n^{k_0+1}[$ we get for
$t\in[t_n^{k_0},t_n^{k_0+1}[$
\begin{equation}\label{eq:cons_en_disc2}
\begin{split}
E_{n,m}(t,u_{n,m}(t))&=E_{n,m}(0,\varphi)\\
&-\sum_{k=1}^{k_0}\int_{t_n^{k-1}}^{t_n^k}\Re\<(i\nabla_x-A(t_n^{k-1}))u_{n,m}(t_n^{k-1}),\partial_tA(t,x)u_{n,m}(t_n^{k-1})\>dt\\
&+O(\frac t n\|\varphi\|^2_{L^2}).
\end{split}
\end{equation}
With this equation we can show that the sequence
$(u_{n,m})_{(n,m)\in\N\times\N}$ is bounded in $H^1_{mg}$. The proof
is a discretization of the proof of Lemma 3.6 in \cite{CaWe88}. Let
$M=2\|\varphi\|_{H^1_{mg}}$ and let $T_{n,m}>0$ the maximal time
such that $1+2m_AT_{n,m}+m_A^2T_{n,m}^2\leq\frac 5 4$ and for
$t\in[0,T_{n,m}[$,
$$\|u_{n,m}\|_{H^1_{mg}}\leq M.$$
Thanks to Propositions \ref{prop:equivnorme},
\ref{prop:estim_nonlin} and Lemma \ref{lem:injecsobmagn} there
exists $K(M)>0$ independent of $n,m\in\N$, such that
$$\|\partial_tu_{n,m}\|_{H^{-1}_{mg}}\leq K(M),\;\forall n,m\in\N,\forall t\in[0,T_{n,m}[ $$
and consequently,
\begin{equation} \label{eq:estL2}
\|u_{n,m}-\varphi\|_{L^2}\leq 2MK(M)t,\;\forall t\in[0,T_{n,m}[.
\end{equation}
On the other hand,
it follows from (\ref{eq:cons_en_disc2}) that
\begin{equation}\label{eq:estim_en_disc}
\begin{split}
\frac 1 2\|(i\nabla_x-&A_n(t))u_{n,m}(t)\|_{L^2}^2\leq \frac 1
2\|(i\nabla_x-A(0))\varphi\|_{L^2}^2
+G_m(u_{n,m})-G_m(\varphi)\\
&-\sum_{k=1}^{n}\int_{t_n^{k-1}}^{t_n^k}\Re\<(i\nabla_x-A(t_n^{k}))u_{n,m}(t_n^{k}),\partial_sA(s,x)u_{n,m}(t_n^{k-1})\>ds\\
&+O(\frac t n\|\varphi\|^2_{L^2}).
\end{split}
\end{equation}
As $\partial_tA$ is bounded, the fourth term of the right hand side
of (\ref{eq:estim_en_disc}) is bounded by $CtM^2$. Moreover it
follows from Lemma \ref{lem:estim_energnl} and estimate
(\ref{eq:estL2}) that
$$|G_m(u_{n,m})-G_m(\varphi)|\leq C(M)(t+t^\nu).$$
Combining these equations with Proposition \ref{prop:equivnorme} we
get
$$\|u_{n,m}\|_{H^1_{mg}}^2\leq
\frac{25}{16}\frac{M^2}2+C(M)(T_{n,m}+T_{n,m}^\nu).$$
 Taking $0<T_{n,m}< T$ with $T$ sufficiently small independently on $n,m$, this
 proves that
\begin{equation}\label{eq:estimapriori}
\|u_{n,m}(t)\|_{L^\infty([0,T],H^1_{mg})}\leq M,\;\forall n,m\in\N
\end{equation}
 Let now $p,q\in\N$, then
 $$i\partial_t(u_{p,m}-u_{q,m})(t)=H_p(u_{p,m}-u_{q,m})(t)+R_{p,q,m}(t)+g_m(u_{p,m}(t))-g_m(u_{q,m}(t))$$
 and $(u_{p,m}-u_{q,m})_{|t=0}=0$, where
 \begin{equation*}
 \begin{split}
 R_{p,q,m}(t)=((A_q-A_p)(i\nabla-A(0))&+(i\nabla-A(0))(A_q-A_p)(t)\\
 &+(A_p^2-A_q^2)(t)+2A(0)(A_q-A_p)(t))u_{q,m}(t).
 \end{split}
 \end{equation*}
 Thanks to Theorem \ref{th:strichmagn}, we have for
 $\widetilde T\in]0,T[$,
 $r\in[,\frac{2n}{n-2}[$ and $\frac 2 q=n(\frac 1 2-\frac 1 r)$,
\begin{equation*}
\begin{split}
\|u_{p,m}&-u_{q,m}\|_{L^q([0,\widetilde T],L^r(\R^n))}\leq
\|R_{p,q,m}\|_{L^{\infty}([0,\widetilde T],L^{2}(\R^n))}\\
&+ C(M)(T
\|u-v\|_{L^{\infty}([0,\widetilde T],L^{2}(\R^n))}+\widetilde T^{\frac{\gamma_2-1}{\gamma_2}}\|u-v\|_{L^{\gamma_2}([0,\widetilde T],L^{r_2}(\R^n))}).
\end{split}
\end{equation*}
On the other hand, $\epsilon>0$ being fixed, for $p,q$ large enough
we have
$$\sup_{(t,x)\in\R\times\R^n}|A_p-A_q|\leq\epsilon.$$
Hence,
$$\|R_{p,q,m}\|_{L^{\infty}([0,\widetilde T],L^{2}(\R^n))}\|\leq
2\epsilon\|u_{q,m}\|_{H^1_{mg}}+C\epsilon\|u_{q,m}\|_{L^2}\leq C
M\epsilon,$$ and for $p,q$ large enough we get
\begin{equation*}
\begin{split}
\|u_{p,m}-u_{q,m}\|_{L^q([0,\widetilde T],L^r(\R^n))}\leq \epsilon&+ C(M)\widetilde T
\|u-v\|_{L^{\infty}([0,\widetilde T],L^{2}(\R^n))}\\
&+C(M)\widetilde T^{\frac{\gamma_2-1}{\gamma_2}}\|u-v\|_{L^{\gamma_2}([0,\widetilde T],L^{r_2}(\R^n))}.
\end{split}
\end{equation*}
This estimate is available, both for $(q,r)=(\infty,2)$ and
$(q,r)=(\gamma_2,\rho_2)$. Summing the two inequalities obtained and
making $\widetilde T>0$ small enough, we get
$$\|u_{p,m}-u_{q,m}\|_{L^q([0,\widetilde T],L^r(\R^n))}\leq 2\epsilon.$$
Therefore, the sequence $(u_{n,m})_{n\in\N}$ converges, as $n$ goes
to infinity, to a limit $u_m\in L^2$ which is solution of
(\ref{eq:nlsmagn_app}). Moreover, as $(u_{n,m})_{n\in\N}$ is bounded
in $H^1_{mg}$ we can suppose that it converges weakly to $u_m$ in
$H^1_{mg}$.

Now let's go back to equation (\ref{eq:cons_en_disc2}). Using the
fact that $u_{n,m}$ converges in $L^2$ and converges weakly in
$H^1_{mg}$ it is no hard to see that
$E_{n,m}(t,u_{n,m})-E_{n,m}(0,\varphi)$ converges as $n\rightarrow\infty$, to
$\Re\int_0^t\<\partial_sA(s)u_m(s),(i\nabla_x-A(s))u_m(s)\>ds$. From
Proposition \ref{prop:estim_nonlin} and weak lower semicontinuity of the
magnetic Sobolev norm $\|(i\nabla-A(t))\,.\|_{L^2}$ it follows that
$$E_m(t,u_m)\leq
E_m(0,\varphi)-\Re\int_0^t\<\partial_sA(s)u(s),(i\nabla_x-A(s))u(s)\>ds.$$

Finally,  $t>0$ being fixed, consider $v_{n,m}(s)=u_{n,m}(t-s)$, which is solution of 
$$i\partial_s v_{n,m}=-\operatorname{H_{A(t-s)}}v_{n,m}-g_m(v_{n,m})$$
with initial data $v_{n,m}(s=0)=u_{n,m}(t)$. Then we can do the same computations as above to get the converse inequality
and hence (\ref{eq:estimenergie_app})  is proved. \fin

\subsubsection{Convergence to the initial problem} In this section,
we show that the sequence $u_m$ converges to a solution of (\ref{eq:nlsmagn}) when $m$ goes to infinity..

\begin{lemma}\label{lem:cauchyL2}
There exists $\tT_2>0$ depending only on $\|\varphi\|_{H^1_{mg}}$ such
that $(u_m)_{m\in\N}$ is a Cauchy sequence in $C([0,\tT_2],L^2)$.
\end{lemma}
\dem The proof is the same as in \cite{CaWe88}, making use of
Theorem \ref{th:strichmagn}, Lemma \ref{lem:estim_energnl},
Proposition \ref{prop:estim_nonlin} and Proposition
\ref{prop:exist_app}. \fin

Now, we can complete the proof of Theorem \ref{th:pbcauchy_nlsmagn}.
This is the same as in \cite{CaWe88} and we recall it for reader
convenience. We denote $u$ the limit of $u_m$ in $C([0,\tT_2],L^2)$.
From estimate (\ref{eq:estimapriori}), it follows that $u\in
L^{\infty}([0,\tT_2],H^1_{mg})$ and by Lemma \ref{lem:injecsobmagn},
$u_m$ converges to $u$ in $C([0,\tT_2],L^r)$ for all $r\geq
2n/(n-2)$. Hence, it follows from Proposition
\ref{prop:estim_nonlin} that $f_m(u_m)$ converges to $f(u)$ in
$C([0,\tT_2],H^{-1}_{mg})$ and $u$ solves (\ref{eq:nlsmagn}) in
$L^{\infty}([0,\tT_2],H^{-1}_{mg})$. Moreover, combining Lemma
\ref{lem:estim_energnl} and Proposition \ref{prop:exist_app} we
prove that
$$E(t,u)=E(0,\varphi)-\Re\int_0^t\<\partial_sA(s)u(s),(i\nabla_x-A(s))u(s)\>ds.$$
This shows that $u\in C([0,\tT_2[,H^1_{mg})$ and hence $u\in
C^1([0,\tT_2[,H^{-1}_{mg})$.

\section{WKB approximation}
In this section we justify WKB approximation for solution of
(\ref{eq:nlsmagn}) when the strength of the magnetic field $b$ goes
to infinity and obtain instability results. We stress our attention
on the case where the magnetic field and the non-linearity have the
same strength; that is we consider the case $\gamma=2$ and search
approximate solution for
\begin{equation}\label{eq:nlswkb}
\left\{\begin{array}{cc}
i\partial_su=\operatorname{H_{A(s)}}u+b^2ug(|u|^2))&\\
u_{|s=0}=a_0(x)e^{ibS(x)}\phantom{******}&
\end{array}
\right.
\end{equation}
where $g$ does not depend on $x$. Remark that with the previous notations, $f=ug(|u|^2)$. In thi section we still ask $f$ to satisfies Assumption \ref{hyp:nonlinearite} and we require additionnaly
\begin{assumption}\label{hyp:nonlinearitebis}
$g\in C^{\infty}(\R_+,\R)$ with $g'>0$.
\end{assumption}

Remark that if we suppose that  $a_0\in H^1$ and $\nabla S+A(0)\in L^2$ then the initial data satisfies
$\|a_0(x)e^{ibS(x)}\|_{H^1_{mg}}=O(b)$. Therefore, under Assumptions
\ref{hyp:dec_magnpot}, \ref{hyp:nonlinearite} and \ref{hyp:nonlinearitebis} it follows from Theorem
\ref{th:pbcauchy_nlsmagn} that there exists a unique solution of
(\ref{eq:nlswkb}) in $C(-T_b,T^b],H^1_{mg})$ with $T_b,T^b\geq C
b^{-\delta}, \delta=\max(2,\frac{2}\alpha)$. In fact this solution
takes a particular form.

\begin{theorem}\label{th:wkbnlsmagn}Let $\sigma>\frac n 2 +2$ and
suppose that Assumptions \ref{hyp:dec_magnpot}, \ref{hyp:nonlinearite} and
\ref{hyp:nonlinearitebis} are satisfied. Assume additionally that
$\partial_tA$ belongs to $H^{\sigma-1}(\R^n)$ for all $t\in\R$ and
take $a_0$ in $H^{\sigma}(\R^n)$ and $S$ such that  $\nabla
S+A(t=0)$ belongs to $H^{\sigma-1}(\R^n)$. Then, there exists $T>0$
and $\alpha_b,\phi_b$ in $C([0,T[,H^{\sigma}(\R^n))\cap
C^1([0,T[,H^{\sigma-1}(\R^n))$ such that
$u(t,x)=\alpha_b(bt,x)e^{ib(S(x)+\phi_b(bt,x))}$ is solution of
(\ref{eq:nlswkb}) on $[0,b^{-1}T]$.
\end{theorem}

\dem We start the proof by a time rescalling leading to a
semiclassical feature. We 
denote $h=b^{-1}>0$ and set $u(s)=v(bs)$. Then equation (\ref{eq:nlswkb}) is equivalent to
\begin{equation}\label{eq:nlswkb_sc}
\left\{\begin{array}{cc}
ih\partial_tv=(ih\nabla_x-A(ht))^2v+vg(|v(t)|^2)&\\
v_{|t=0}=a_0(x)e^{ih^{-1}S(x)}\phantom{**********}&
\end{array}
\right.
\end{equation}

We follow the general method initiated by Grenier \cite{Gr98} for the semiclassical Schr\"odinger
equation and look for a phase and an amplitude depending on the parameter $h$. Putting the ansatz
$v(t,x)=\alpha_h(t,x)e^{ih^{-1}\phi_h(t,x)}$ in the equations (\ref{eq:nlswkb_sc}) we get
\begin{equation}\label{eq:systeiktrans1}
\left\{\begin{array}{cc}
\partial_t\phi_h+|\nabA\phi_h|^2+g(|\alpha_h|)^2=0\phantom{*********}&\\
\partial_t\alpha_h+\nabA\phi_h.\nabla\alpha_h+div(\nabA\phi_h)\alpha_h=ih\Delta\alpha_h&
\end{array}
\right.
\end{equation}
where $\nabA\phi=(\nabla_x\phi+A(ht))$. Next we set $\varphi_h(t,x)=\nabA\phi_h(t,x)\in\R^n$ and differentiate the above eikonal equation with respect to $x$. We obtain
\begin{equation}\label{eq:systeiktrans2}
\left\{\begin{array}{cc}
\partial_t\varphi_h+2\varphi_h.\nabla \varphi_h +2g'(|\alpha_h|^2)\Re(\overline{\alpha_h}\nabla\alpha_h)=h\partial_tA(ht,x)\phantom{}&\\
\partial_t\alpha_h+\varphi_h.\nabla\alpha_h+div(\varphi_h)\alpha_h=ih\Delta\alpha_h\phantom{*********}&
\end{array}
\right.
\end{equation}
 Separating real and imaginary parts of $\alpha_h=\alpha_{1,h}+i\alpha_{2,h}$,
(\ref{eq:systeiktrans2}) becomes
\begin{equation}\label{eq:systhyp}
\partial_tw_h+\sum_{j=1}^nA_j(w_h)\partial_{x_j}w_h=hLw_h+\nu_h
\end{equation}
with
\begin{equation}\label{eq:defvect}
w_h=\left(\begin{array}{c}
\alpha_{1,h}\\
\alpha_{2,h}\\
\varphi_{1,h}\\
\vdots\\
\varphi_{n,h}
\end{array}
\right)
,
\nu_h=\left(\begin{array}{c}
0\\
0\\
h\partial_tA_1(ht,x)\\
\vdots\\
h\partial_tA_n(ht,x)
\end{array}
\right)
\end{equation}
\begin{equation}
L=\left(\begin{array}{ccccc}
0&-\Delta&0\\
\Delta&0&0\\
0&0&0_{n\times n}
\end{array}
\right)
\end{equation}
and
\begin{equation}
A_j(w)=\left(\begin{array}{ccccc}
\varphi_{j,h}&0&\alpha_1&\ldots&\alpha_1\\
0&\varphi_{j,h}&\alpha_2&\ldots&\alpha_2\\
2g'\alpha_1&2g'\alpha_2&v_j&0&0\\
\vdots&\vdots&0&\ddots&0\\
2g'\alpha_1&2g'\alpha_2&0&0&\varphi_{j,h}
\end{array}
\right)
\end{equation}
This system has the same form as in \cite{Gr98}, \cite{Ca06} with
the exception of the source term $\nu_h$ in right hand side of
(\ref{eq:systhyp}) and the initial data. Thanks to the assumptions,
$\nu_h$ belongs to $H^{\sigma-1}(R^n)$, whereas the initial condition
in (\ref{eq:nlswkb_sc}) yields
\begin{equation}\label{eq:condinit_hyp}
w_h(t=0)=\left(\begin{array}{c}
\Re a_o\\
\Im a_0\\
\partial_{x_1}S+A_1(0)\\
\vdots\\
\partial_{x_n}S+A_n(0)
\end{array}
\right)
\end{equation}
which belongs to $H^{\sigma-1}(\R^n)$.

On the other hand, thanks to the assumption on $g'$, the system (\ref{eq:systhyp}) can be
symmetrized by
\begin{equation}
S=\left(\begin{array}{cc}
I_2&0\\
0&\frac 1 {g'}I_n
\end{array}
\right)
\end{equation}
which is symmetric and positive. It follows from general theory of
hyperbolic systems that the problem (\ref{eq:systhyp}) together with
initial condition (\ref{eq:condinit_hyp}) has a unique solution
$w_h\in L^{\infty}([0,T_h],H^{\sigma-1})$ for some $T_h>0$.

Hence, we have to bound $T_h$ from below by a constant independent
of $h$. This is done by computing classical energies estimates as in
\cite{Gr98}, \cite{Ca06}, and using the fact that $\partial_t A$ as
well as $\nabla_x S+A(0)$ belong to $H^{\sigma-1}.$

Finally we define $\alpha_h$ and $\phi_h$ by
$\alpha_h=w_{1,h}+iw_{2,h}$ and $$\phi_h=S(x)-\int_0^t
|\varphi_{h}|^2+f(|\alpha_h|^2)ds.$$ By construction, $\phi_h$ belongs to
$L^2$. Moreover, a simple calculus shows that $\nabla_x\phi_h=\varphi_h-A(ht)$
belongs to $H^{\sigma-1}$ so that  $\phi_h$ is in fact in
$H^{\sigma}$. Going back to the equation on $\alpha_h$ and making
energies estimates we show that $\alpha_h\in H^{\sigma}$. 
Finally, it a direct calculus shows that $(\alpha_h,\phi_h)$ defined above solves (\ref{eq:systeiktrans1})

\fin

\begin{remark}
The above solution belongs to the magnetic sobolev space $H^1_{mg}$.
Indeed,
$$(i\nabla_x-bA)(\alpha_b
e^{ib\phi_b})=(i\nabla\alpha_b-b(\nabla\phi_b+A)\alpha_b)e^{ib\phi_b}$$
belongs to $L^2$. Therefore the solution built in Theorem \ref{th:wkbnlsmagn}
coincide with the one of Theorem \ref{th:pbcauchy_nlsmagn}.
\end{remark}

With Theorem \ref{th:wkbnlsmagn} in hand it is easy to prove
instability results.
\begin{proposition}\label{prop:instab}
Let $\sigma>\frac n 2 +2$ and let $A$ satisfy the assumptions of Theorem \ref{th:wkbnlsmagn}. Suppose that $S$ is such that $\nabla
S+A(t=0)$ belongs to $H^{\sigma-1}(\R^n)$. Then, there exists $a_0$ and $\widetilde a_{0,b}$ in $H^{\sigma}(\R^n)$ and $0<t_b<Cb^{-1}$ such that 
$$\|a_0-\widetilde a_{0,b}\|_{L^2}\rightarrow 0\text{ as }b\rightarrow\infty$$
and the solutions $u_b$ (resp. $\widetilde u_b$) associated to (\ref{eq:nlswkb}) with initial data $a_0e^{ibS(x)}$ (resp. $\widetilde a_0e^{ibS(x)}$) satisfy
$$\|u_b-\widetilde u_b\|_{L^{\infty}([0,t_b],L^2)}\geq 1.$$
\end{proposition}
\dem
It is a straightforward consequence of Theorem \ref{th:wkbnlsmagn} and the methods of \cite{Ca06}.
\fin
\bibliographystyle{amsplain}
\bibliography{ref}

\end{document}